\documentclass[11pt]{amsart}
\usepackage{graphicx,calc,amscd}
\usepackage{epsfig,psfig}
\usepackage{float}
\usepackage{labelfig}
\usepackage{amsmath}
\usepackage{amsfonts}
\usepackage{amssymb}

\newcommand{\R}{{\mathbb R}}

\newcommand{\GG}{{\mathbb G}}

\newcommand{\BB}{{\mathcal B}}

\newcommand{\T}{{\mathcal T}}
\newcommand{\M}{{\mathcal M}}
\newcommand{\PP}{{\mathcal P}}
\newcommand{\Y}{{\mathcal Y}}

\newcommand{\fix}{\mathrm{Fix}}

\renewcommand{\tt}{\mbox{\rm Tr}_{2^d|2}\,}

\newtheorem{theorem}{Theorem}[section]
\newtheorem{lemma}[theorem]{Lemma}
\newtheorem{corollary}[theorem]{Corollary}
\newtheorem{proposition}[theorem]{Proposition}
\newtheorem{remark}[theorem]{Remark}

\newtheorem{definition}[theorem]{Definition}
\begingroup\makeatletter\ifx\SetFigFont\undefined%
\gdef\SetFigFont#1#2#3#4#5{%
 \reset@font\fontsize{#1}{#2pt}%
 \fontfamily{#3}\fontseries{#4}\fontshape{#5}%
 \selectfont}%
\fi\makeatother\endgroup%
\newenvironment{example}[1][Example]{\par\normalfont
  \trivlist\item[\hskip\labelsep\itshape #1:]\ignorespaces
  }

\begin{document}
\title{A geometric characterization of orientation reversing involutions}
\author{Antonio F. Costa and Hugo Parlier}
\address{Departamento de Matematicas Fundamentales\\
Facultad de Ciencias\\
Universidad Nacional de Educaci\'{o}n a Distancia\\
Madrid 28040\\
SPAIN} \email{acosta@mat.uned.es} \email{hugo.parlier@epfl.ch}

\thanks{The first author was supported in part by BFM 2002-04801, the second author was supported by SNFS grant number PBEL2-106180.}
\date{\today}
\begin{abstract}We give a geometric characterization of compact Riemann surfaces
admitting orientation reversing involutions with fixed points.
Such surfaces are generally called real surfaces and can be
represented by real algebraic curves with non-empty real part. We
show that there is a family of disjoint simple closed geodesics
that intersect all geodesics of a partition at least twice in
uniquely right angles if and only if the involution exists. This
implies that a surface is real if and only if there is a pants
decomposition of the surface with all Fenchel-Nielsen twist
parameters equal to $0$ or $\frac{1}{2}$.
\end{abstract}
\maketitle

\section{Introduction}\label{Sect:S1}

A smooth complex projective algebraic curve $C$ can be represented
by a compact Riemann surface $S$ (i.e. an orientable compact
surface with a conformal structure). The curve $C$ can be
described by real polynomial equations if and only if the surface
$S$ admits an orientation reversing involution $\sigma$ defined by
complex conjugation. If the set of real points of $C$ is not empty
then the fixed point set of $\sigma$ is non-empty or equivalently
the field of meromorphic functions of $(S,\sigma)$ is real (i. e.
in such field $-1$ is not a sum of squares). Hence a surface $S$
admitting an orientation reversing involution $\sigma$ with
$\fix(\sigma)\neq\emptyset$ is called a \textit{real Riemann
surface} and $\sigma$ is a real form on $S$.

In this article we deal with one of the fundamental problems in
this subject: to decide whether a complex algebraic curve admits a
real form, i. e. whether a complex algebraic curve can be
described using polynomials with real coeficients. Equivalently
the problem is to describe in $\mathcal{M}_{g} $
(moduli space of surfaces of genus $g$) the real moduli $\mathcal{M}%
_{g}^{\mathbb{R}}$ whose points are the real Riemann surfaces. The
main result determines, in the Teichm\"{u}ller space of Riemann
surfaces of genus $g,$ $\T_{g}$, a subspace $\R\T_{g}$ that
projects on $\mathcal{M}_{g}^{\mathbb{R}}$. We use the classic
Fenchel-Nielsen parametrization of $\T_{g}$, namely the parameters
are the collection of lengths of geodesics in a pants
decomposition and twist parameters which describe how the pairs of
pants are pasted together. Hence, the parametrization is not
homogeneous in the nature of the parameters. We show that all
points of $\mathcal{M}_{g}^{\mathbb{R}}$ can be represented by
elements of $\T_g$ where all the twist parameters are zero or
$1/2$ (and vice-versa), hence the Fenchel-Nielsen parametrization
restricted to $\R\T_g$ produces a parametrization where all
parameters are of homogeneous nature. Real moduli has recently
been studied by several authors, see for instance \cite{busi01},
\cite{de04}, and \cite{nabook04}, with new and classical
applications.

In order to prove the above result we obtain a geometric
characterization of the Riemann surfaces admitting an orientation
reversing involution by showing that such surfaces are exactly the
Riemann surfaces having a pants decomposition with with zero or
$1/2$ twist parameters. The pants decomposition turns out to be
not only invariant by the orientation reversing involution but
furthermore the involution induces the identity of the graph of
such a decomposition. This essential fact is no longer possible if
one considers an orientation reversing involution without fixed
points, as the number of disjoint simple closed geodesics that are
left invariant by such an involution cannot exceed $g+1$. We also
mention that the usage of invariant pants decompositions for real
Riemann surfaces can be found in Buser and Sepp\"{a}la's paper
\cite{buse92}, but they allow the involution to define a
non-trivial combinatorial involution on the graphs of such
decompositions.

The main result is Theorem 3.2 where we prove that a Riemann
surface is real if and only if there exists a set $\BB$ of
disjoint simple closed geodesics that intersect all geodesics of a
pants decomposition at least twice at perpendicular angles. There
are similar geometrical characterizations for surfaces admitting
conformal involutions, obtained by Maskit \cite{mas00} and
Schmutz-Schaller (\cite{sc001}, \cite{sc002}). Using
uniformization groups, there are other characterizations of real
Riemann surfaces, see for instance \cite{sibn65}, \cite{bujsi85},
and \cite{he87}.

Finally, in Section 4, we obtain as a consequence information on
the upper bound of the distance between fixed points of an
orientation reversing involution and other points of the surface.
We also provide an example illustrating this last result.

\section{Preliminaries}
Our main object of study is a compact Riemann surface $S$  (i.e.
an orientable compact surface with a conformal structure) of genus
$g\geq 2$. It can be endowed with a hyperbolic metric which we
shall denote $d(\cdot,\cdot)$. Generally, curves and geodesics
will be considered primitive and non-oriented, and will be seen as
point sets on $S$. Occasionally, we will need to consider surfaces
with boundary, and the signature of such a surface will be denoted
$(g,k)$ where $g$ is the topological genus and $k$ is the number
of simple boundary curves. Unless specified, boundary curves will
be considered geodesic. The following propositions concern well
known properties of curves on hyperbolic surfaces (i.e.
\cite{bubook}, pp. 19-23).

\begin{proposition}\label{prop:simplecurves}
Let $S$ be a hyperbolic surface. Let $\alpha,\beta$ be disjoint
simple closed geodesics on $S$. Let $c$ be a simple path from
$\alpha$ to $\beta$. Then in the free homotopy class of $c$ with
endpoints gliding on $\alpha$ and $\beta$, there exists a unique
shortest curve, denoted $\GG(c)$, which meets $\alpha$ and $\beta$
perpendicularly. Furthermore, if $\tilde{c}$ is also a simple path
from $\alpha$ to $\beta$ such that $c\cap \tilde{c}=\emptyset$,
then either $\GG(c)=\GG(\tilde{c})$ or
$\GG(c)\cap\GG(\tilde{c})=\emptyset$.
\end{proposition}

In the case of simple closed curves, the corresponding proposition
is the following. Note that we call a piecewise geodesic boundary
curve $c$ {\it convex} if for all points $p$ of $c$ the interior
angle $\theta_p$ verifies $\theta_p\leq \pi$.

\begin{proposition}\label{prop:simpleclosedcurves}
Let $S$ be a compact hyperbolic surface (possibly with piecewise
geodesic convex boundary) and let $c$ be a homotopically
non-trivial simple closed curve on $S$. Then $c$ is freely
homotopic to a unique simple closed geodesic, denoted $\GG(c)$.
The curve $\GG(c)$ is either contained in $\partial S$ or $\GG(c)
\cap
\partial S = \emptyset$. If $c$ is a non-smooth boundary
component, then $\GG(c)$ and $c$ bound an embedded annulus (see
figure \ref{fig:2}).
\end{proposition}

A collection $\PP$ of $3g-3$ disjoint simple closed geodesics is
called a partition or the geodesics of a pants decomposition, and
$S\setminus \PP$ is a collection of $2g-2$ surfaces of signature
$(0,3)$, commonly called pairs of pants or $Y$-pieces. Following
proposition \ref{prop:simplecurves}, between two distinct boundary
geodesics of a $Y$-piece $\Y$, there is a unique geodesic path
perpendicular to the boundary geodesics. The three perpendicular
paths defined in this way are disjoint, and by
cutting along them one obtains two anticonformal isometric hexagons.\\

There are different ways to define what are commonly called
{\it{twist parameters}}. The definition we will use is the
following, which is relative to a given partition.

\begin{definition}
Let $\gamma$ be a geodesic in $\PP$ with a given orientation. Let
$\Y_1=(\alpha,\beta,\gamma)$ and $\Y_2=(\alpha',\beta',\gamma)$ be
the two $Y$-pieces pasted along $\gamma$. Let
$\ell_{\alpha\gamma}$ be the perpendicular between $\alpha$ and
$\gamma$ on $\Y_1$. Let $\ell_{\alpha'\gamma}$ be the
perpendicular between $\alpha'$ and $\gamma$ on $\Y_2$. Let $p$ be
the intersection point between $\gamma$ and $\ell_{\alpha\gamma}$
and let $q$ be the intersection point between $\gamma$ and
$\ell_{\alpha'\gamma}$. Let $\ell$ be the length of the path on
$\gamma$ from $p$ to $q$ following the orientation of $\gamma$.
The twist parameter $t_\gamma$ along $\gamma$ is defined as the
quantity $\frac{\ell}{\gamma}$.
\end{definition}

The above definition fixes twist parameters in the interval
$[0,1[$, and the definition can be extended to a parameter that
lies in $\R$. Using the extended version of twist parameters and
the lengths of geodesics in a partition, one obtains the
Fenchel-Nielsen parametrization of the space of marked Riemann
surfaces of genus $g$ (\cite{feni48}, \cite{fenibook}, or
Teichm\"uller space of genus $g$ surfaces, denoted $\T_g$. This
shows that $\T_g$ is homeomorphic to $(\R^{+*})^{3g-3}\times
\R^{3g-3}$. The Fenchel-Nielsen parametrization of Teichm\"uller
space is relative to a graph describing the partition where the
graph associated to a pants decomposition has vertices which
correspond to $Y$-pieces and two vertices are connected if and
only if there is a common boundary for the two corresponding
$Y$-pieces. Restricting the twist parameters to the interval
$[0,1[$ still ensures us that we have at least one representative
(in fact an infinity) of each conformal equivalence class of
surfaces of genus $g$. The space of all conformal equivalence
classes of surfaces of genus $g$ is called the Moduli
space and is denoted $\M_g$.\\

An {\it involution} is an isometry of the surface onto itself that
is of order $2$. An involution can either be orientation reversing
or not. The following remarks hold for any orientation reversing
involution, although in the sequel we will see that there are
fundamental differences between
the case where an involution is with or without fixed points.\\

The following proposition is an extension of what is generally
called Harnack's theorem (\cite{wephd}) and can be found in
\cite{krma82}.

\begin{proposition}\label{prop:harnackgen}
If a surface $S$ admits $\sigma$, an orientation reversing
involution, then the fixed point set of $\sigma$ is a set of $n$
disjoint simple closed geodesics
${\BB}=\{\beta_1,\hdots,\beta_n\}$ with $n\leq g+1$. In the case
where the set $\BB$ is separating, then $S\setminus \BB$ consists
of two connected components $S_1$ and $S_2$ such that $\partial
S_1 =\partial S_2 = \BB$ and $S_2=\sigma(S_1)$. If not, then $\BB$
can be completed by either a set $\alpha$ which consists of one or
two simple closed geodesics such that $\BB \cup \alpha$ has the
properties described above (with the important difference that
$\alpha$ does not contain any fixed points of $\sigma$). Each of
the simple closed geodesics in $\alpha$ is globally fixed by
$\sigma$.
\end{proposition}

The following proposition concerns any simple closed geodesic
fixed by an orientation reversing involution.

\begin{proposition}\label{prop:fixedgeod}
Let $S$ be a surface admitting an orientation reversing involution
$\sigma$. Let $\gamma \subset S$ be a simple closed geodesic such
that $\sigma(\gamma)=\gamma$. If $\gamma$ does not contain any
fixed points of $\sigma$, then for all $p\in \gamma$ the image
$\sigma(p)$ of $p\in \gamma$ is the point on $\gamma$
diametrically opposite from $p$.
\end{proposition}

\begin{proof} For any $p\in \gamma$, if the points $p$ and $\sigma(p)$
are distinct then they separate $\gamma$ into two geodesic arcs.
The image of one of the arcs is either the other arc or is
globally fixed. Since $\sigma$ is an isometry, the mid point of
the fixed arc must be fixed by $\sigma$ and the result follows.
\end{proof}

\section{Orientation reversing involutions with fixed points}

This section will be devoted to a geometric characterization of
surfaces admitting an orientation reversing involution with fixed
points. We denote by $\sigma$ an orientation reversing involution
with fixed points. Thus the fixed point set of $\sigma$ is a
non-empty set of $n$ disjoint simple closed geodesics
${\BB}=\{\beta_1,\hdots,\beta_n\}$ with $n\leq g+1$. By cutting
$S$ along $\BB$, one obtains either one or two surfaces with
boundary. The following lemma concerns surfaces with boundary, and
will be useful in the sequel.

\begin{lemma}\label{lemma:cut} Let $S$ be a Riemann surface with non-empty boundary with boundary geodesics
$\beta_1,\hdots,\beta_k$. Let $c_1,\hdots,c_j$ be disjoint simple
geodesic paths perpendicular to the boundary of $S$. Then $j\leq
6g-6+3n$ and the set of $c_i$s can be completed by
$c_{l+1},\hdots,c_{6g-6+3n}$ such that the new set verifies the
same conditions as above. \end{lemma}

\begin{proof}
The idea of the proof is to show that $c_1,\hdots,c_j$ can be
completed into a set that decomposes the surface into hexagons.
Once this is obtained, we will show that all further simple
geodesic paths perpendicular to boundary must intersect curves in
this set, and that the number of simple paths is exactly $6g-6+3n$.\\

The connected components of $S\setminus \{c_1,\hdots,c_j\}$ are
orientable surfaces with boundary. All boundary curves are
piecewise geodesic and are composed of an even number of simple
geodesic paths. Furthermore, an orientation can be given to the
boundary curve such that the oriented angles of intersection are
always $\frac{\pi}{2}$. The boundary geodesics are either trivial
(they are the boundary of a geodesic polygon) or are freely
homotopic to a simple closed geodesic on $S$ that is
disjoint with the considered boundary.\\

First let us deal with the case when a boundary curve, say $b$, is
the boundary of a hyperbolic polygon. The polygon is a convex
right-angled $2i$-gon with $i\geq 3$.\\

\begin{figure}[h]
\leavevmode \SetLabels
\L(.4*.93) $b_1$\\
\L(.57*.025) $b_7$\\
\L(.38*.08) $b_9$\\
\L(.67*.47) $b_5$\\
\L(.33*.4) $$\\
\L(.65*.4) $$\\
\endSetLabels
\begin{center}
\AffixLabels{\centerline{\epsfig{file = 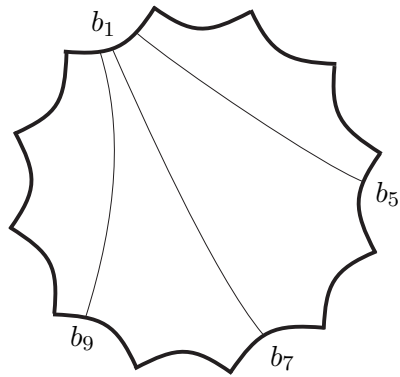,
width=5cm,angle= 0}}}
\end{center}
\caption{The decomposition of a $12$-gon} \label{fig:twelvegon}
\end{figure}

If $b_1, \hdots, b_{2i}$ are the simple geodesic paths that
compose $b$, then notice that the numbering can be chosen such
that for $j$ odd, $b_j$ is an element of $\beta$, and for $j$
even, $b_j$ is a $c_{j_0}$ for some $j_0$. From $b_1$, consider
the unique geodesic perpendicular paths to each $b_j$ for $5\leq
j\leq 2i-3$ odd. (Put each one of these paths in our collection of
$c_j$s.) By cutting along these paths, one obtains $i-2$
hyperbolic
right-angled hexagons as desired.\\

\begin{figure}[h]
\leavevmode \SetLabels
\L(.4*.7) $\GG(b)$\\
\L(.33*.4) $$\\
\L(.65*.4) $$\\
\endSetLabels
\begin{center}
\AffixLabels{\centerline{\epsfig{file = 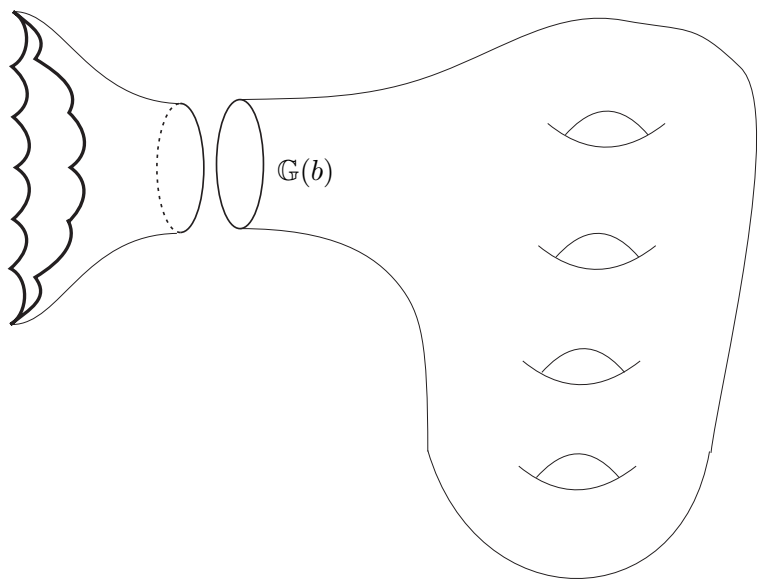,
width=10cm,angle= 0}}}
\end{center}
\caption{The cylinder seen on the surface} \label{fig:2}
\end{figure}

If $b$ is not a trivial boundary curve, and $\GG(b)$ the unique
simple closed geodesic on $S$ that is freely homotopic to it, then
$b$ and $\GG(b)$ are the boundary curves of a topological cylinder
embedded in $S$.

\begin{figure}[h]
\leavevmode \SetLabels
\L(.46*.89) $$\\
\L(.575*.08) $$\\
\L(.46*.08) $$\\
\L(.575*.89) $$\\
\L(.33*.4) $$\\
\L(.65*.4) $$\\
\endSetLabels
\begin{center}
\AffixLabels{\centerline{\epsfig{file = 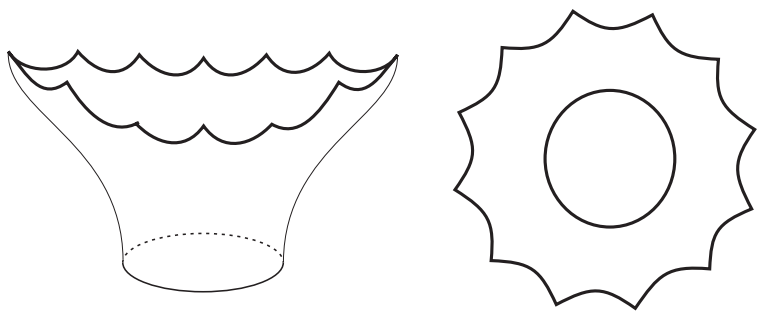,
width=10cm,angle= 0}}}
\end{center}
\caption{Two ways of seeing the cylinder enclosed by $b$ and
$\GG(b)$} \label{fig:geodesicboundary}
\end{figure}

Say $b_1, \hdots, b_{2i}$ are the simple geodesic paths that
compose $b$. As before, the numbering can be chosen such that for
$j$ odd, $b_j$ is an element of $\beta$, and for $j$ even, $b_j$
is a $c_{j_0}$ for some $j_0$.\\

The next step is to cut hyperbolic hexagons out of the cylinders
in order to reduce the boundary of the cylinder to a curve with
two geodesic components. Consider the unique geodesic
perpendicular path $c$ from $\GG(b)$ to $b_1$. By cutting along
$c$ and $\GG(b)$ one obtains a $(2i+4)$-gon, and following the
same procedure as for polygons, one obtains $i$ hexagons. Now
pasting along the path $c$ we have a cylinder with two boundary
components. One of these is a geodesic $2$-gon. Notice that in the
$2$-gon, one of the geodesic paths is an embedded
element of $\beta$, and the other is a boundary to boundary path.\\

\begin{figure}[h]
\leavevmode \SetLabels
\L(.28*.78) $b_1$\\
\L(.315*.68) $c$\\
\L(.46*.08) $$\\
\L(.575*.89) $$\\
\L(.33*.4) $$\\
\L(.65*.4) $$\\
\endSetLabels
\begin{center}
\AffixLabels{\centerline{\epsfig{file = 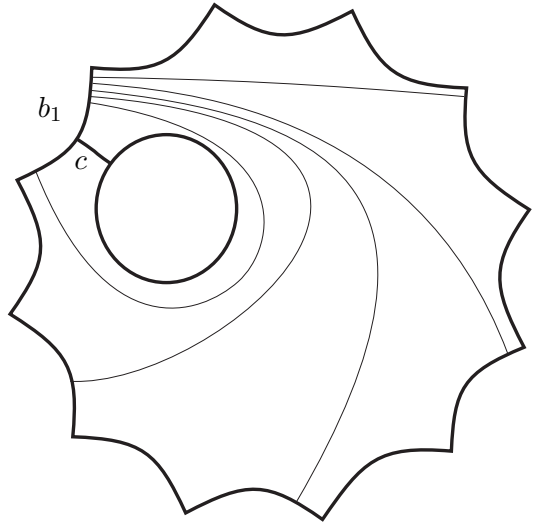,
width=7cm,angle= 0}}}
\end{center}
\caption{The procedure to obtain a piecewise geodesic boundary
with two components} \label{fig:twelvegonboundary}
\end{figure}

For all $b$, geodesic boundary paths, consider $\GG(b)$. Cutting
along these, one obtains a collection of surfaces with simple
geodesic boundaries, and cylinders obtained as above. On a surface
with simple geodesic boundaries, consider a perpendicular geodesic
path $c$ between two boundary geodesics, say $\gamma$ and
$\delta$. If $\gamma$ and $\delta$ are connected components of
$\beta$, then the path is admissible. If not, at least one (say
$\gamma$) is perpendicular to the boundary geodesic of a cylinder
as described above. On this cylinder consider the unique
perpendicular geodesic path $\tilde{c}$ from the embedded element
of $\beta$ to $\gamma$. If $c\cap \gamma=p$ and $\tilde{c}\cap
\gamma=q$, then consider the (or a) shortest path on $\gamma$
between $p$ and $q$, say $\bar{c}$. Then $\GG(c\cup\bar{c}\cup
\tilde{c})$ is a simple perpendicular path as desired. Notice that
the procedure works even if $\gamma=\delta$ is a geodesic embedded
in a cylinder. After each cut, the boundary curves may have to be
reduced to geodesic $2$-gon boundary as explained above, or a
geodesic polygon may have been
cut off.\\

This procedure can be pursued until the surface has been cut into
hexagons. We shall call the completed set of paths along which $S$
has been cut $\{c_i\}_{i=1}^{\tilde{n}}$. Notice that on such a
hexagon, if a side is an embedded part of $\beta$, then it's
adjacent sides are two (not necessarily distinct) $c_i$s (and
vice-versa). Although there are still perpendicular paths between
sides on hexagons, they are necessarily between opposite edges,
thus from embedded parts of $\beta$ to a $c_i$. This fact proves
the maximality of the set $\{c_i\}$. We must now count the
hexagons. As the hyperbolic area of a right-angled hexagon is
$\pi$, and the area of $S$ is $2\pi(2g-2+n)$, it follows that
there are $2(2g-2+n)$ hexagons. Each hexagon has three $c_i$s and
each path is counted twice. It follows that $\tilde{n}=6g-6+3n$
and the lemma is proved.
\end{proof}

We can now proceed to the main result.

\begin{theorem}\label{thm:FP}
Let $S$ be a Riemann surface of genus $g$. Then $S$ admits an
orientation reversing involution $\sigma$ with fixed points if and
only if there exists a set $\BB$ of disjoint simple closed
geodesics and a partition $\PP$ such that all intersections
between the geodesics of two sets are perpendicular and such that
all geodesics in $\PP$ intersect $\BB$ at least twice.
Furthermore, $\fix(\sigma)=\beta$ and $\sigma(\gamma)=\gamma$ for
all $\gamma\in \PP$.
\end{theorem}

\begin{proof}
Let us begin by proving the existence of $\BB$ and $\PP$ as above
if $S$ is real. Denote by $\sigma$ the real involution of $S$ and
by $\BB=\{\beta_1,\hdots,\beta_n\}$ the fixed point set of
$\sigma$. There are two cases to consider. The first case is when
$\BB$ separates $S$ into two surfaces $S_1$ and $S_2$ such that
$\sigma(S_1)=S_2$. Notice that $g$ and $n$ are of different parity
both $S_1$ and $S_2$ are of signature $(g',n)$ where
$g'=(g-n+1)/2$. We can apply lemma $\ref{lemma:cut}$ and on $S_1$
there exists a set $c_1,\hdots,c_{3g-3}$ of disjoint simple
geodesic perpendicular to boundary paths. As all points of
$\BB=\partial S_1$ are fixed by $\sigma$, it follows that for all
$i\in \{1,\hdots,3g-3\}$, $\gamma_i=c_i\cup \sigma(c_i)$ is a
simple closed geodesic. The set
$\PP=\{\gamma_1,\hdots,\gamma_{3g-3}\}$ is a partition, as all
$\gamma_i$s are disjoint. Each $\gamma_i$ intersects the set $\BB$
twice and in right angles.

The idea for the case where $\BB$ is not separating is identical,
but necessitates a prelude. Consider the set $\alpha$ defined as
in proposition \ref{prop:harnackgen}, which consists of one or two
disjoint simple closed geodesics such that $\sigma(\alpha)=\alpha$
and $S\setminus \BB \setminus \alpha$ has two connected components
$S_1$ and $S_2$ such that $\sigma(S_1)=S_2$. On $S_1$, consider a
perpendicular simple geodesic path $c$ oriented from a geodesic in
$\BB$, say $\beta_1$, to a geodesic in $\alpha$, say $\alpha_1$.
If $\alpha_1$ is given an orientation, the piecewise geodesic path
$c \alpha_1 c^{-1}$ is freely homotopic to exactly one purely
simple geodesic perpendicular path from $\beta_1$ to $\beta_1$,
say $c_1=\GG(c\alpha_1 c^{-1})$. The path $c_1$ is separating for
$S_1$, cuts $\beta_1$ in two, and the surface cut off by $c_1$ is
of signature $(0,2)$ with boundaries $\alpha_1$ and a piecewise
geodesic curve consisting of part of $\beta_1$, and $c_1$. The set
$\gamma_1=c_1\cup \sigma(c_1)$ is a separating simple closed
geodesic on $S$, and cuts off a $Q$-piece containing $\alpha_1$.
If $\alpha$ contains another simple closed geodesic, say
$\alpha_2$, then the same procedure can be applied to cut off
another $Q$-piece containing $\alpha_2$ with boundary geodesic
$\gamma_2$, which intersects another element of $\BB$ in two
points in a perpendicular fashion, and verifies
$\sigma(\gamma_2)=\gamma_2$.

\begin{figure}[h]
\leavevmode \SetLabels
\L(.19*.615) $\beta_1$\\
\L(.25*.8) $c_1$\\
\L(.23*.11) $\sigma(c_1)$\\
\L(.375*.615) $\alpha_1$\\
\L(.33*.4) $$\\
\L(.65*.4) $$\\
\endSetLabels
\begin{center}
\AffixLabels{\centerline{\epsfig{file = 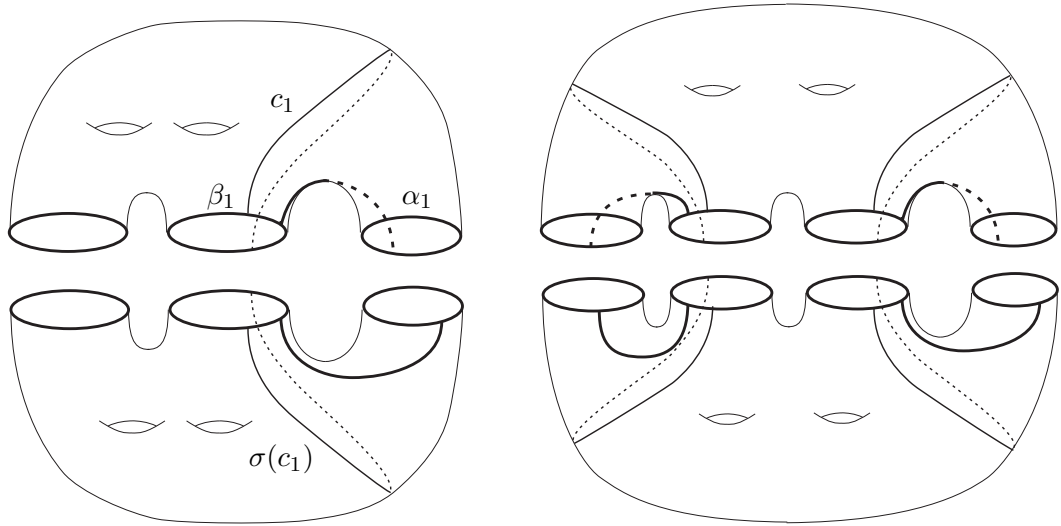,
width=14cm,angle= 0}}}
\end{center}
\caption{The two possible cases if $\beta$ is non-separating}
\label{fig:betanonseparating}
\end{figure}

The geodesic $\gamma_1$ (resp. the geodesics $\gamma_1$ and
$\gamma_2$) separate $S$ into a $Q$-piece $Q_1$, and a surface
$\tilde{S}$ of signature $(g-1,1)$ (resp. into two $Q$-pieces
$Q_1$, $Q_2$, and a surface $\tilde{S}$ of signature $(g-2,2)$).
The proof of when $\BB$ is separating applies to $\tilde{S}$, as
the set $\BB\mid_{\tilde{S}}$ is now separating.  We shall proceed
to deal with $Q_1$ (and $Q_2$ if necessary).\\

Consider $Q_1$ (the proof for $Q_2$ is identical if necessary).
The geodesic $\gamma_1$ is exactly the boundary of $Q_1$, and by
construction, $\beta_1\mid_{Q_1}=h$ is a non-trivial simple
geodesic perpendicular path from $\gamma_1$ to $\gamma_1$. Cutting
along $h$, one obtains a surface $C$ of signature $(0,2)$, with
piecewise geodesic boundary consisting of $c_1\cup h_1$ and
$\sigma(c_1)\cup h_2$ where $h_1$ and $h_2$ are the two copies of
$h$ on $C$. Notice that $h_1=\sigma(h_2)$ and that $\alpha_1$ is
an interior separating geodesic of $C$.\\

\begin{figure}[h]
\leavevmode \SetLabels
\L(.128*.92) $h_1$\\
\L(.565*.92) $h_1$\\
\L(.16*.63) $l_1$\\
\L(.38*.3) $l_2$\\
\L(.28*.65) $C_1$\\
\L(.28*.25) $C_2$\\
\L(.39*.03) $h_2$\\
\L(.83*.03) $h_2$\\
\L(.27*.375) $\alpha_{11}$\\
\L(.27*.525) $\alpha_{12}$\\
\L(.69*.65) $\GG(\tilde{c_1})$\\
\L(.66*.28) $\GG(\tilde{c_2})$\\
\endSetLabels
\begin{center}
\AffixLabels{\centerline{\epsfig{file =
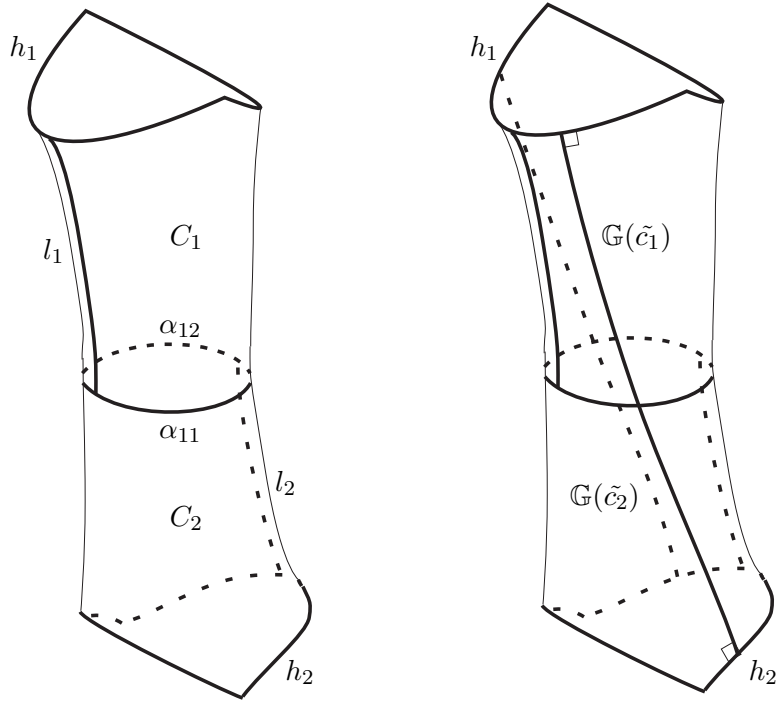, width=10cm,angle= 0}}}
\end{center}
\caption{$C$ and certain curves}
\label{fig:fixedpointproofcylinder}
\end{figure}

Denote by $C_1$ (resp. $C_2$) the connected component of
$C\setminus \alpha_1$ with boundary $c_1\cup h_1$ (resp. the
connected component of $C\setminus \alpha_1$ with boundary
$\sigma(c_1)\cup h_2$). Notice that $\sigma(C_1)=C_2$. The
shortest path $l_1$ on $C_1$ between $\alpha_1$ and $h_1$ is the
only simple geodesic perpendicular path to both $\alpha_1$ and
$h_1$. Denote by $l_2$ the corresponding path on $C_2$ (between
$h_2$ and $\alpha_1$). Because both $l_1$ and $l_2$ are unique on
$C_1$ and $C_2$, it follows that $\sigma(l_1)=l_2$. Denote by
$p_1$ and $p_2$ the points $l_1\cap \alpha_1$ and $l_2\cap
\alpha_1$. Cutting along $\alpha$ from $p_1$ to $p_2$, one obtains
two equal length geodesic arcs (as $p_2 = \sigma(p_1)$) which we
shall denote by $\alpha_{11}$ and $\alpha_{12}$. Now consider the
two homotopically distinct non-oriented paths $l_1 \alpha_{11}
l_2$ (resp. $l_1 \alpha_{12} l_2$) and denote them by
$\tilde{c_1}$ and $\tilde{c_2}$. It is easy to see that the
geodesics in their free homotopy class are disjoint, simple,
perpendicular to both $h_1$ and $h_2$ and that
$\tilde{c_2}=\sigma(\tilde{c_1})$. On $S$ they are both
perpendicular paths from $\beta$ to $\beta$ and it follows that on
$\tilde{c_1}\cup \tilde{c_2}$ is a simple closed geodesic that
intersects $\beta_1$ twice at right angles, and does not further
intersect any element of $\BB$.\\

\begin{figure}[h]
\leavevmode \SetLabels
\L(.46*.89) $$\\
\L(.575*.08) $$\\
\L(.46*.08) $$\\
\L(.575*.89) $$\\
\L(.33*.4) $$\\
\L(.65*.4) $$\\
\endSetLabels
\begin{center}
\AffixLabels{\centerline{\epsfig{file = 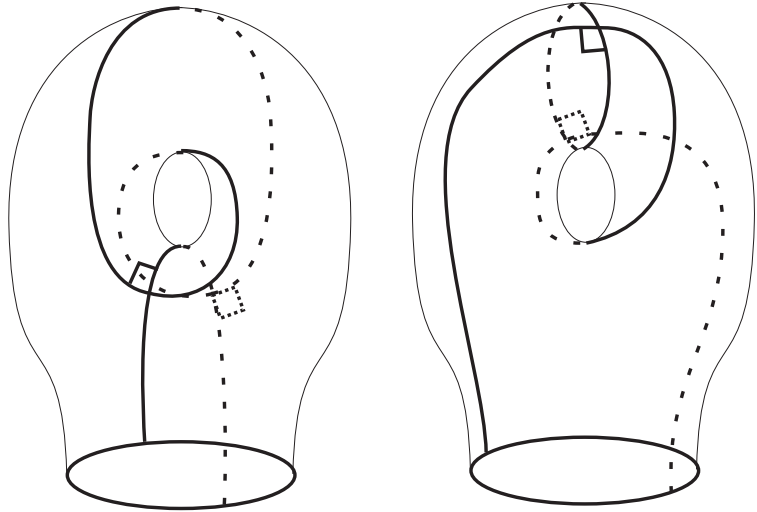,
width=10cm,angle= 0}}}
\end{center}
\caption{The partition geodesic obtained seen two different ways}
\label{fig:fixedpointproofcylinder}
\end{figure}

We shall now prove the reciprocal. Consider a set $\beta$ and a
partition $\PP$ that intersects $\beta$ as in the hypotheses. As
each geodesic in $\PP$ intersects $\beta$ at least twice, then for
a pair of pants $\Y$ in the underlying pants decomposition, the
three perpendicular paths from boundary to boundary are subsets of
$\beta$. To prove this, consider the boundary geodesics of $\Y$,
say $\gamma_1,\gamma_2,\gamma_3$. The conditions impose that
$\beta\mid_\Y$ must consist of simple perpendicular paths between
the boundary geodesics. Thus, the connected components of
$\beta\mid_\Y$ are either one of the three perpendiculars
mentioned above (type $1$), or simple perpendicular paths that
whose end points lie on a same boundary geodesic (type $2$). If
$\beta\mid_\Y$ consists only of paths of type $1$, then as each
boundary geodesic intersects at least two of them, the three
perpendicular paths of type $1$ are contained in $\beta\mid_\Y$.
If $\beta\mid_\Y$ contains an element of type $2$, then suppose
that it intersects $\gamma_1$. In this case, then $\gamma_2$ and
$\gamma_3$ can only intersect one perpendicular path (both of type
$1$) and this is a contradiction.\\

\begin{figure}[h]
\leavevmode \SetLabels
\L(.46*.89) $$\\
\L(.575*.08) $$\\
\L(.46*.08) $$\\
\L(.575*.89) $$\\
\L(.33*.4) $$\\
\L(.65*.4) $$\\
\endSetLabels
\begin{center}
\AffixLabels{\centerline{\epsfig{file = 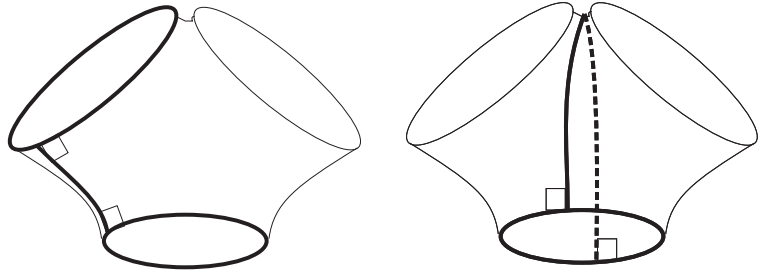,
width=10cm,angle= 0}}}
\end{center}
\caption{Perpendicular paths of type $1$ and $2$}
\label{fig:betanonseparating}
\end{figure}

It follows that every pair of pants $\Y_i \in S\setminus \PP$ is
divided into two anticonformal isometric hexagons, say $H_i$ and
$\tilde{H}_i$, by $\beta \mid_{\Y_i}$. For every $\Y_i \in \PP$,
consider the local orientation reversing involution $\sigma_i$
that takes a point on $H_i$ to its corresponding point
$\tilde{H}_i$. The involutions $\sigma_i$ can be extended to act
on $\PP$. This gives a unique involution $\sigma_\PP$ that acts on
the geodesics of $\PP$. From the set $\sigma_i$s and $\sigma_\PP$,
we obtain an application $\sigma$ on $S$ defined by
$\sigma(p)=\sigma_i(p)$ if $p\in \Y_i$. It is straightforward to
see that $\sigma$ is an orientation reversing involution that
verifies $\fix(\sigma)=\beta$ and $\sigma(\gamma)=\gamma$ for all
$\gamma \in \PP$.
\end{proof}

\begin{figure}[h]
\leavevmode \SetLabels
\L(.46*.89) $$\\
\L(.575*.08) $$\\
\L(.46*.08) $$\\
\L(.575*.89) $$\\
\L(.33*.4) $$\\
\L(.65*.4) $$\\
\endSetLabels
\begin{center}
\AffixLabels{\centerline{\epsfig{file = 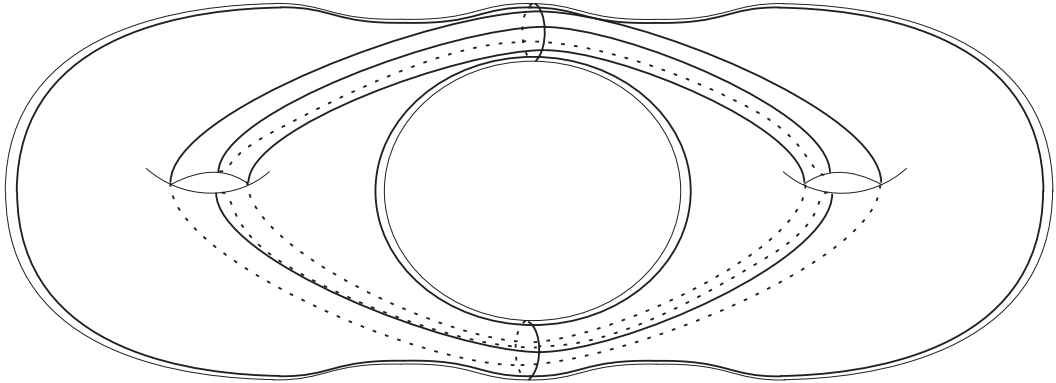,
width=14cm,angle= 0}}}
\end{center}
\caption{A genus $3$ surface admitting an involution with two
fixed geodesics} \label{fig:genus3partition}
\end{figure}

This result can be seen in function of Fenchel-Nielsen parameters.

\begin{corollary}\label{cor:FPFN}
Let $S$ be a Riemann surface of genus $g$. Then $S$ admits an
orientation reversing involution $\sigma$ with fixed points if and
only if there exists a pants decomposition such that the
Fenchel-Nielsen twist parameters are all equal to $0$ or
$\frac{1}{2}$.
\end{corollary}

\begin{proof}
If such an involution exists, then consider the partition $\PP$
and set of geodesics $\BB$ guarantied by the previous theorem.
Consider the set of $Y$-pieces $\{\Y_i\}$ that form $S\setminus
\PP$. For all $i\in 1,\hdots,2g-2$, the set $\BB \mid_{\Y_i}$ is
the set of the perpendicular paths of type $1$ on $\Y_i$. Hence
the only possibility is that the twist parameters lie in the set
$\{0,\frac{1}{2}\}$ if $\BB$ is to consist of simple closed
geodesics. Reciprocally, if there exists a pants decomposition
with twist parameters $0$ and $\frac{1}{2}$, then the
perpendicular simple paths between distinct boundary geodesics of
the pairs of pants, form a collection of simple closed geodesics
which are disjoint and intersect all geodesics in the underlying
partition exactly twice.
\end{proof}

\begin{remark}
The involution acts as the identity on the graph induced by the
pants decomposition described above.
\end{remark}

\section{Further consequences of the geometric characterization}

One of the main consequences of theorem \ref{thm:FP} is that it
gives a very precise image of surfaces which allow such
involutions. This vision for instance allows the following
proposition which concerns the distance between fixed points of an
orientation reversing involution and other points of the surface.

\begin{proposition}\label{prop:disk}
The following inequality is always true:
\begin{equation}\label{eq:disk1}
\max_{p\in S}d_S(\fix(\sigma),p)>\frac{\ln 3}{2}.
\end{equation}
Reciprocally, for any $\varepsilon>0$ and any $g\geq 2$, there
exists a surface $S_\varepsilon$ of genus $g$ with an orientation
reversing involution with fixed points (of any given type) such
that
\begin{equation}\label{eq:disk2}
\max_{p\in
S_\varepsilon}d_{S_\varepsilon}(\fix(\sigma_\varepsilon),p)<\frac{\ln
3}{2} + \varepsilon.
\end{equation}

\end{proposition}

\begin{proof}
In \cite{pa042}, the existence of a disk of radius $\frac{\ln
3}{2}$ on the complement of a partition is shown. As the geodesics
composing the fixed point set of $\sigma$ can be completed into a
partition, it follows that on $S\setminus \fix(\sigma)$, there is
a disk of radius $\frac{\ln 3}{2}$ and equation \ref{eq:disk1}
follows. For equation \ref{eq:disk2}, it suffices to construct an
example. For a given topological type of involution $\sigma$,
theorem \ref{thm:FP} ensures us of the existence of a partition
$\PP$ such that all perpendicular paths between boundary geodesics
of the underlying $Y$-pieces are elements of $\fix(\sigma)$. We
shall consider the lengths of the geodesics in $\PP$ as free
parameters to describe the surface $S_\epsilon$ without modifying
the twist parameters. Using the methods in \cite{pa042}, for any
$\varepsilon>0$, there exists a constant $C_\varepsilon$ such that
if one chooses all geodesics in $\PP$ of length shorter than
$C_\epsilon$, one obtains a surface $S_\varepsilon$ without any
disks of radius $\frac{\ln 3}{2}$ embedded in
$S_\varepsilon\setminus\PP\setminus \fix(\sigma)$. This proves
equation \ref{eq:disk2}. For more clarity see the example that
follows the proof.

\end{proof}

\begin{example}\label{ex:sharp} In order to illustrate proposition \ref{prop:disk},
consider the following example.

\begin{figure}[H]
\begin{center}
\epsfig{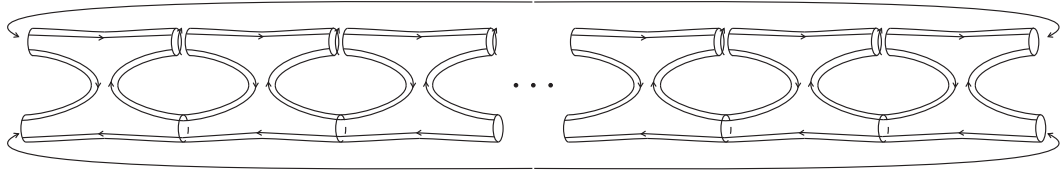}
\end{center}
\caption{An example in genus $g$}\label{fig:finalexample}
\end{figure}

Notice that the twist parameters of certain geodesics have been
left free. If, for instance, all twist parameters are considered
to be $\frac{1}{2}$, then we have an example where $\fix(\sigma)$
is a single simple closed geodesic. If all twist parameters are
equal to $0$, then $\fix(\sigma)$ consists of $g+1$ simple closed
geodesics. The twist parameters can be chosen in order to create
an example with any topological type of $\sigma$.
\end{example}

The fixed point geodesics of $\sigma$ can be chosen as short as
wanted. Hence, the collar theorem (i.e. \cite{bu78},\cite{kee74})
implies that for any constant $C$, there exists a surface $S_C$
admitting an orientation reversing involution $\sigma_C$ with
fixed points such that
\begin{equation*}
\max_{p\in S_C}d_{S_C}(\fix(\sigma_C),p)>C.
\end{equation*}

In corollary \ref{cor:FPFN}, we show the existence of a partition
$\PP$ such that each simple closed geodesic is left invariant by
$\sigma$. In \cite{buse92}, one of the main results concerns
partitions of surfaces with an orientation reversing involution.
It is shown that a partition $\tilde{\PP}$ can be chosen such that
$\sigma(\tilde{\PP})=\tilde{\PP}$, and that $\max_{\gamma \in
\tilde{\PP}}\ell(\gamma) \leq 21 g$ where $g$ is the genus. The
existence of such a constant for arbitrary Riemann surfaces of
same genus was originally proved in \cite{be85}. However, although
the set $\tilde{\PP}$ is globally fixed by $\sigma$, each geodesic
in $\tilde{\PP}$ is not necessarily fixed. The graph of the
underlying pants decomposition is invariant, but the involution
does not necessarily act as the identity on it. Using the collar
theorem as above, it is easy to see that we can not find a bound
on the lengths of the geodesics in $\PP$ in function of the genus
$g$, as was done in \cite{buse92}. If however, one imposes a lower
bound on the lengths of geodesics in $\fix (\sigma)$, then using
the methods exposed in \cite{buse92}, for given genus, a partition
$\PP$ that verifies the conditions of corollary \ref{cor:FPFN},
can be chosen with bounded length.\\

\bibliographystyle{plain}
\def\cprime{$'$}

\end{document}